\documentclass[a4paper,11pt]{amsart}

\usepackage{fullpage,amsxtra,amssymb,amsthm,amsmath}
\usepackage[latin1]{inputenc}

\newcommand{\notation}[2]{\quad {#1}& : &{#2}\\}

\renewcommand{\leq}{\leqslant}
\renewcommand{\geq}{\geqslant}

\newcommand{\sheaf}[1]{\mathcal{{#1}}}

\def\stacksum#1#2{{\stackrel{{\scriptstyle #1}}
{{\scriptstyle #2}}}}

\def\triplesum{\mathop{\sum\sum\sum}\limits}

\newcommand{\Cc}{\mathbf{C}}
\newcommand{\Aa}{\mathbf{A}}

\newcommand{\Zz}{\mathbf{Z}}

\newcommand{\Rr}{\mathbf{R}}

\newcommand{\Qq}{\mathbf{Q}}
\newcommand{\Fp}{\mathbf{F}}

\newcommand{\barre}[1]{\overline{{#1}}}

\newcommand{\mods}[1]{\,(\mathrm{mod}\,{#1})}

\newcommand{\ra}{\rightarrow}

\DeclareMathOperator{\spec}{Spec}

\DeclareMathOperator{\frob}{Fr}
\DeclareMathOperator{\Gal}{Gal}

\DeclareMathOperator{\Tr}{Tr}

\newcommand{\eps}{\varepsilon}

\newcommand{\demi}{{\textstyle{\frac{1}{2}}}}

\DeclareMathSymbol{\gena}{\mathord}{letters}{"3C}
\DeclareMathSymbol{\genb}{\mathord}{letters}{"3E}

\theoremstyle{plain}
\newtheorem{theorem}{Theorem}

\newtheorem{lemma}[theorem]{Lemma}
\newtheorem{corollary}[theorem]{Corollary}

\newtheorem{proposition}[theorem]{Proposition}

\theoremstyle{remark}
\newtheorem{remark}[theorem]{Remark}

\theoremstyle{definition}

\newtheorem{example}[theorem]{Example}

\newcommand{\tformula}[1]{\text{``{{#1}}''}}

\begin{document}

\title{Exponential sums over definable subsets of finite fields}
\author{E. Kowalski}
\address{Universit\'e Bordeaux I - A2X\\
351, cours de la Lib\'eration\\
33405 Talence Cedex\\
France}
\email{emmanuel.kowalski@math.u-bordeaux1.fr}
\begin{abstract}
We prove some general estimates for exponential sums over subsets of
finite fields which are definable in the language of rings. This
generalizes both the classical exponential sum estimates over
varieties over finite fields due to Weil, Deligne and others, and the
result of Chatzidakis, van den Dries and Macintyre
concerning the number of points of those definable sets.
As a first application, there is no
formula in the language of rings that defines for infinitely many
primes an ``interval'' in $\Zz/p\Zz$ that is neither bounded nor with
bounded complement. 
\end{abstract}

\maketitle

\section{Introduction}

Exponential sums are ubiquitous in analytic number theory, in various
shape and forms. A basic type is a sum
\begin{equation}\label{eq-gen-exp}
S_f(M,N)=\sum_{M\leq n<M+N}{e(f(n))},
\end{equation}
where $e(z)=e^{2i\pi z}$ and $f$ is some real-valued function. These
tend to arise naturally in 
any asymptotic counting problem, as ways to express the secondary
terms after isolating a ``main term'' and the basic goal is to establish
some form of cancellation, of the type
\begin{equation}\label{eq-gen-bound}
\sum_{M\leq n<M+N}{e(f(n))}\ll N\theta(N)^{-1},
\end{equation}
where the saving $\theta(N)$ from the trivial bound $N$ is a positive
increasing function with $\theta(N)\ra +\infty$ as $N\ra
+\infty$. 
Evidently, it must be the case that $f$ varies ``fast enough'' for
such an estimate to hold.
\par
Various highly ingenious methods have been developed to deal with the
distinct possible types of phase functions $f$; the names of Weyl, van
der Corput and Vinogradov in particular are attached to the most
classical ideas (see e.g.~\cite[\S 8]{ant}).
It was however discovered 
that this type of analytic questions could
sometimes be attacked using  highly involved algebraic tools:
if the interval of summation is of the type $0\leq n<p$, where $p$ is
prime, and if $f(n)=g(n)/p$, where 
$g$ is a polynomial or a rational function, the best general results
come from an interpretation as an exponential sum over the finite
field $\Zz/p\Zz$. 
\par
Indeed, one introduces the ``companion'' sums
$$
S_{\nu}=\sum_{x\in \Fp_{p^{\nu}}}{e\Bigl(\frac{\Tr f(x))}{p}\Bigr)},
$$
for $\nu\geq 1$, where $\Fp_{p^{\nu}}$ is a field with $p^{\nu}$ elements,
$\Tr\,:\, \Fp_{p^{\nu}}\ra \Fp_p=\Zz/p\Zz$ being the trace map. Although
$S_{\nu}$, $\nu\geq 2$, never (?) has any interpretation in analytic number
theory, it is the properties of the generating function
$$
Z(T)=\exp\Bigl(\sum_{\nu\geq 1}{\frac{S_{\nu}}{\nu}T^{\nu}}\Bigr)
$$
which are fundamental in understanding the original sums. In this
context, this was first recognized and developed by 
A. Weil, who proved for instance that for a fixed (non-constant)
function $g\in \Zz[X]$ one has 
$$
S_{\nu}(p)\ll p^{\nu/2}
$$
for all primes $p$ and $\nu\geq 1$ (with possibly few well-understood
exceptions), with an implied constant depending only on $g$. See
e.g~\cite[\S 11]{ant} for a description of the elementary approach of
Stepanov and~\cite[\S 11.11]{ant} for a first survey of the more
advanced cohomological methods of Grothendieck, Deligne, Katz and
others. 
\par
In terms of applications to analytic number theory, it is clear that
the potential of the more advanced results has not yet been fully exploited;
there are a number of reasons for this, not only the complexity of the
algebraic geometry involved (although that is certainly a factor), but
also the difficulty of bringing a natural problem to a position where
the Riemann Hypothesis for varieties over finite fields can be applied
successfully: the reader need only look at the proof of the Burgess
estimate for short character sums (see e.g.~\cite[12.4]{ant}) to see
what ingenuity may be required; also the comments
in~\cite[11.12]{ant} explain how the question of uniformity in
parameters and ``flexibility'' in the shape of the sums can be crucial
matters. 
\par
In this paper, we describe a new general estimate for exponential
sums over finite fields which combines quite efficiently the
cohomological methods (as ``black-box'') and some results and
techniques of logic to give estimates where
the summation set in the finite field is much more general than the
algebraic sets that are usually considered. We hope that this added
flexibility will make it suitable for applications to analytic number
theory; also the statement is, in itself, quite elementary with very
few conditions, and this may also make it appealing to readers
without a great experience in algebraic geometry (at the price of
learning, or remembering, a few notions of logic...)
\par
See Section~\ref{sec-statement} for the general statement, which is
preceded by a section recalling the precise formulas permitted as
summation conditions. As a sample result, for the very important case
of a sum in one variable, one gets the following: 

\begin{theorem}\label{th-one-var}
Let $\varphi(x)$ be a first-order formula in the language
$(0,1,+,-,\cdot)$ of rings.\footnote{\ We recall the precise definition
  in Section~\ref{sec-definable}.} For every ring $A$, let 
$$
\varphi(A)=\{x\in A\,\mid\, \varphi(x)\text{ holds}\}.
$$
\par
Let $f$, $g\in \Qq(X)$ be rational functions with $f$
non-constant. Let $N\geq 1$ be the product of primes $p$ such that $f$
modulo $p$ is constant. Then there exists a constant $C\geq 0$, depending only
on $\varphi$ and the degree of the numerator and denominator of $f$
and $g$ such that for any prime $p$ and any multiplicative character
$\chi$ modulo $p$ we have 
\begin{equation}\label{eq-one-var}
\Bigl|\sum_{\stacksum{x\in \varphi(\Zz/p\Zz)}
{f(x), g(x)\text{ defined}}}{\chi(g(x))
e\Bigl(\frac{f(x)}{p}\Bigr)}\Bigr|
\leq C(p,N)^{1/2}\sqrt{p}.
\end{equation}
\end{theorem}

Compared to the classical sums above, the point is that the summation condition
can be quite complicated, involving arbitrary entanglements of
quantifiers (in first-order predicates, i.e., applied to elements of
the field). 
One may also wonder if in fact the bound is really non-trivial (what if the
number of points is usually of size $p^{1/4}$, for instance?), but in
fact, as proved in~\cite{zoe} and as we will explain again in detail
below, the number of points of summation is either $\leq A$ or
or $\geq cp$, for some $A\geq 1$ and $c>0$ depending only on the
formula $\varphi$. And one
should keep in mind that if this were applied to a problem of analytic
number theory, whether this is efficient or not would most often be
obvious from the final result anyway. 
\par
\medskip
\textbf{Notation.}  As already mentioned, we  denote
$e(z)=e^{2\pi iz}$. We denote by $\Fp_q$ a finite field with $q$
elements and usually $p$ is its characteristic. For a finite set $X$,
$|X|$ denotes its cardinality. By $f\ll g$ for $x\in X$, or $f=O(g)$
for $x\in X$, where $X$ is an arbitrary set on which $f$ is defined, we mean
synonymously that there exists a constant $C\geq 0$ such that
$|f(x)|\leq Cg(x)$ for all $x\in X$. The ``implied constant'' is any
admissible value of $C$. It may depend on the set $X$ which is always
specified or clear in context.
For notation and conventions concerning logical formulas, see the
beginning of the next section.
We use elementary scheme-theoretic language for our algebraic
geometry (see e.g.~\cite[II]{hart}); in particular, an algebraic
variety over a field $F$ or over $\Zz$ is simply a separated scheme of
finite type over $F$ or $\Zz$, and in fact only affine schemes will
occur (so separatedness is automatic); so a variety is not necessarily
reduced or irreducible. We write either $V_A$ or $V/A$ to indicate
that a scheme is defined over a ring $A$. (The choice of $A$ is
sometimes important to indicate precise dependency for constants that
occur.) 

\section{Definable sets}
\label{sec-definable}

Since the paper involves a fairly unusual mixture of analytic number
theory and logic (also algebraic geometry, but the latter is kept
essentially inside the proofs), we start by recalling what are
precisely the formulas which define the summation sets we will
consider. Of course logicians can skip this section without loss,
unless they wish to make sure that the author does not speak utter
nonsense. 
\par
A \emph{term} in the language $(0,1,+,-,\cdot)$ of rings is simply a
polynomial $f\in\Zz[x_1,\ldots,x_n]$ with integer coefficients, where
the $x_i$ are variables; an \emph{atomic formula} $\varphi$ is a 
formula of the type $f=g$ where $f$ and $g$ are polynomials
(possibly involving distinct sets of variables). Given an atomic formula
$\varphi$, a ring $A$, and assignments of elements $x_i=a_i\in A$ to the
variables involved, the formula $\varphi(a)$ with $a=(a_i)$
substituted for the variables is \emph{satisfied} in $A$, denoted 
$$
A\models \varphi(a)
$$
if the equality which ``is'' $\varphi$ holds when the variables are
given the values $a_i$. 
\par
Next the (first-order) \emph{formulas} in the language of rings are
built from atomic formulas by induction using the additional logical
symbols $\neg$, $\wedge$, $\vee$, and quantifiers $\exists$,
$\forall$: atomic formulas are formulas by definition and if
$\varphi$ and $\varphi_i$, $i\in I$, are formulas, with $I$ finite,
then  
$$
\neg \varphi_i\quad\quad\quad
\bigwedge_{i\in I}{\varphi_i}\quad\quad\quad
\bigvee_{i\in I}{\varphi_i}
$$
are also formulas, and if $x$ is any variable then
$$
\exists x\ \varphi\quad\text{ and } \quad \forall x\ \varphi
$$
are also formulas. Implication $\ra$ and equivalence $\leftrightarrow$
are defined as abbreviations: 
\begin{gather*}
\varphi_1\ra \varphi_2 \text{ means } (\neg\varphi_1)\vee \varphi_2,\\
\varphi_1\leftrightarrow \varphi_2\text{ means }
(\varphi_1\ra \varphi_2)\wedge (\varphi_2\ra \varphi_1),
\end{gather*}
(and parenthesizing can be introduced for clarity). 
\par
In an obvious way, the relation 
$A\models \varphi(a)$ defined for an atomic formula is extended to any
formula by induction using the usual meanings of the symbols: $\wedge$
as ``and'', $\vee$ as ``or'', $\neg$ as ``not''. The quantifiers are
always extended to elements of $A$ only (not to subsets, not to
elements of other rings than $A$). For instance the torsion subgroup
of invertible elements in $A$ is \emph{not} definable in the
first-order language of rings since the exponent can not be bounded a
priori. 
\par
Given a formula $\varphi(x,y)$, where $x=(x_1,\ldots, x_n)$ and
$y=(y_1,\ldots, y_m)$ are (disjoint) tuples of variables, and given a
ring $A$ and $y\in A^m$, we put
\begin{equation}\label{eq-def-set}
\varphi(A,y)=\{x\in A^n\,\mid\, A\models \varphi(x,y)\},
\end{equation}
in other words, $\varphi(A,y)$ is the set of $n$-tuples in $A$ which
satisfy the formula $\varphi$ for the given value $y$ of the
parameter. Such a set is called a \emph{definable set} (with
parameters); if $m=0$ (no parameters), then the set is also called
$\emptyset$-definable. 
\par
If $\varphi$ is an atomic formula $f=g$, the assignment $A\mapsto
\varphi(A)$ is simply the functor that associates its $A$-valued
points $V(A)$ to the scheme $V=\spec(\Zz[X]/(f-g))$ over $\Zz$ given
by the equation 
$f=g$. So in general definable sets can be seen as a (substantial)
generalization of algebraic varieties. 
\par
For instance, the set of squares in $A$ is $\emptyset$-definable by
the formula $\exists y\ x=y^2$, but the assignment $A\mapsto
\varphi(A)$ is not a functor (e.g. because an element
can become a square after extension to a larger field, which can not
happen to points of schemes). Here is more complicated example:
for $d\geq 1$ an integer, the set of irreducible
polynomials of degree $\leq d$ in $A[X]$, identified with a subset of
$A^{d+1}$ by the coefficients, is $\emptyset$-definable (for every
$j$, $k\geq 1$ with $j+k\leq d$, write existential quantifiers on
coefficients of two polynomials of degree $j$ and $k$ with product
equal to the given one). 
\par
\medskip
See e.g.~\cite[1.3,2.1]{hodges} or~\cite[6.1]{fried-jarden} for more
details (involving more general languages) and more examples. 
\par
A final notation: to define a formula $\varphi(x)$, we will sometimes
write 
$$
\varphi(x)=\tformula{$x$ satisfies such and such property},
$$
(and will either explain or leave to the reader to check that the
property thus stated in informal manner can be written as a
first-order formula in the language of rings), or 
$$
\varphi(x)\,:\, \psi_1(x)\ldots
$$
to indicate that the formula $\varphi(x)$ is defined to be the
expression after ``:'', which will usually  be a combination of
various bits and pieces; for instance
$$
\varphi(x)\,:\, (\exists y\ x=y^2)\vee (\exists y\ x=y^3)\vee 
(\forall y\ \exists z\ y^2+x=z^2).
$$

\section{Exponential sums over definable sets in finite fields}
\label{sec-statement}

This section defines the exponential sums we want to consider. We
generalize from the context described in the introduction by 
introducing formulas with parameters $\varphi(x,y)$ where
$x=(x_1,\ldots, x_n)$, $n\geq 1$, are the variables and
$y=(y_1,\ldots,y_m)$, $m\geq 0$, are the parameters. This  
formula is still assumed to be in the first-order language of
rings. As in the previous section, we let
$$
\varphi(A,y)=\{x\in A^n\,\mid\, A\models \varphi(x,y)\}
$$
for any ring $A$ and parameter $y\in A^m$. 
\par
We consider especially finite fields $\Fp_q$ with $q$ elements. Assume
that for all $q$ in some subset of the powers of primes we have chosen
an additive character
$$
\psi\,:\, \Fp_q\ra \Cc^{\times}
$$
and a multiplicative character
$$
\chi\,:\, \Fp_q^{\times}\ra \Cc^{\times}
$$
(which of course depend on $q$, although it is not indicated
in the  notation). We extend $\chi$ to $\Fp_q$ by putting
$\chi(0)=0$, except when $\chi$ is the trivial character, in which
case $\chi(0)=1$.
\par
Let $f_1$, $f_2$, $g_1$, $g_2\in \Zz[X]$ be polynomials in
$X=(X_1,\ldots, X_n)$, with $f_2$ and $g_2$ non-zero, and let
$f=f_1/f_2$, $g=g_1/g_2$ (as rational functions). We assume that for
all the $q$ under consideration the
formula $\varphi(x,y)$ satisfies
\begin{equation}\label{eq-no-poles}
\Fp_q\models (\varphi(x,y)\ra (f_2(x)g_2(x)\not=0)),
\end{equation}
and if necessary, this can be achieved by replacing $\varphi(x,y)$ by
the formula 
$$
\varphi(x,y)\wedge (\neg f_2(x)=0)\wedge (\neg g_2(x)=0),
$$
which ``restricts'' to the points which are not poles of $f$ and
$g$. This may introduce a further dependency of the results
below on $f_2$ and $g_2$, but that is of course not surprising, and it
will be clear that such dependency is really only in terms of the
degree of $f_2$ and $g_2$.
\par
Now finally we introduce the following general exponential sums over a
definable set:
\begin{equation}\label{eq-exp-def}
S(y,\varphi,\Fp_q)=\sum_{x\in \varphi(\Fp_q,y)}{
\psi(f(x))\chi(g(x))},
\end{equation}
for all $q$ for which the data is defined.
\par
These generalize the more classical exponential sums over the
$\Fp_q$-points of an algebraic variety $V/\Zz$, which corresponds to
the case where the formula $\varphi(x,y)$ is the conjunction of the atomic
formulas which ``are'' the equations of $V$. In that situation we also
denote
\begin{equation}
\label{eq-exp-var}
S(y,V,\Fp_q)=\sum_{x\in V(\Fp_q)}{\psi(f(x))\chi(g(x))}.
\end{equation}
\par
\medskip
The natural goal is to describe a non-trivial upper bound for
$S(y,\varphi,\Fp_q)$ when possible, as explicit as possible in its
dependencies. For applications to analytic number theory, it is
natural to look primarily at the so-called ``horizontal'' case, i.e.,
when $\Fp_q=\Fp_p$ is the prime field, and our statements are
skewed to this case (assuming for instance that $p$ is large enough,
instead of $q$ large enough, for some condition to hold).
\par
Here are two fairly simple examples that follow easily from our
results.

\begin{theorem}\label{th-exp-def}
With data as described above, assume that the additive characters
$\psi$ are non-trivial. There 
exist constants $C\geq 0$ and $\eta>0$, 
depending only on the formula $\varphi(x,y)$ and the degrees of the
polynomials $f_1$, $f_2$, $g_1$, $g_2$ such that for any prime $p$
and any parameter $y\in\Fp_p^m$, we have
$$
\Bigl|\sum_{x\in \varphi(\Fp_p,y)}{\psi(f(x))\chi(g(x))}\Bigr|
\leq Cp^{-1/2}\sum_{x\in\varphi(\Fp_p,y)}{1}
$$
\emph{unless} there exists $c\in \Fp_p$
such that 
$$
|\{x\in\varphi(\Fp_p,y)\,\mid\, f(x)=c\}|\geq \eta 
|\varphi(\Fp_p,y)|.
$$
\end{theorem}

The result is stated in this manner in order to make quite clear what
the saving is (i.e., about $p^{1/2}$) compared to the trivial
estimate. We will also recall during the course of the proof the main
theorem of~\cite{zoe} which gives the approximate value of the number
of points of summation. The condition for the estimate to hold is
sufficient but not necessary (see below for examples). It states
intuitively that if there is no cancellation, then $f$ must be
constant on a subset of $\varphi(\Fp_p,y)$ which has ``positive
density''. 
\par
This condition may seem difficult to check but in applications it
should be the case that the ``degenerate'' cases are fairly
obvious, and can be dealt with separately. The second result takes a
different approach and is a ``baby'' version of the stratification
results of Katz-Laumon (compare with~\cite{fouvry-katz}).

\begin{theorem}\label{th-twists}
Let $\varphi(x)$ be a formula without parameters with $x=(x_1,\ldots,
x_n)$.
Define the additive characters $\psi$ by $\psi(x)=e(\Tr x/p)$, 
and let $\chi$, $f_i$, $g_i$ be as above. There exist constants
$C\geq 0$, $D\geq 0$, depending only on $\varphi(x)$ and, in the case
of $C$, on the degrees of $f_i$, $g_i$ such that for all primes $p$,
all $h=(h_1,\ldots, h_n)\in \Fp_p^n$ except at most $Dp^{n-1}$
exceptions we have
$$
\Bigl|\sum_{x\in\varphi(\Fp_p)}{\chi(g(x))
e\Bigl(\frac{f(x)+h_1x_1+\cdots +h_nx_n}{p}\Bigr)}\Bigr|
\leq C\Bigl(1+p^{-1/2}\sum_{x\in\varphi(\Fp_p)}{1}\Bigr).
$$
\end{theorem}

(Intuitively,  a sum may not exhibit cancellation, but most of its 
``twists'' will; adding $1$ on the right-hand side takes care of those
$p$ where $|\varphi(\Fp_p)|$ is small).
\par
These theorems will be derived from more precise ``structural'' results
concerning the exponential sum and its companions, which are of
independent interest and may be useful for deeper studies
in some cases. 
\par
\medskip
Before going into details, a very important remark is that we have
\emph{not} found a new 
source of cancellation in exponential sums: the saving will come by
application of the Riemann Hypothesis for varieties over finite fields
(Deligne's theorem). However, the point is that we have a very general
result, with great flexibility, and with intrinsic uniformity in
parameters, so Theorem~\ref{th-exp-def} has the potential of being a very
efficient ``black-boxing'' of Deligne's result. Certainly in tricky
cases it might not be easy to perform the reduction to varieties over
finite fields explicitly.

\begin{remark}
We have not stated results with purely multiplicative sums. Although
they can be analyzed by the method of Sections~\ref{sec-decomp}
and~\ref{sec-estimates} this tends to 
reveal trickier aspects which make simple statements for a general
summation set hard to state (in other words, it seems one needs a
deeper understanding of the structure of $\varphi(\Fp_q,y)$). We will
give a few examples in the next section.
\end{remark}

\section{Examples}
\label{sec-comments}

We will now give a few simple examples of the exponential sums we have
in mind, and make some general comments. We make forward references to
the structural results of the next sections, as the examples
illustrate various aspects involved.

\begin{example}
The simplest example of a formula defining a set which is not a
variety is probably the formula
$$
\varphi(x)\,:\, \exists y\ x=y^2
$$
which characterizes the squares in a ring. Thus, if the characteristic
$p$ is odd, we have
$$
|\varphi(\Fp_q)|=\frac{q+1}{2},
$$
and since $x\mapsto x^2$ is two-to-one except at $0$, for any
complex numbers $\beta(x)$ defined for $x\in\Fp_q$ we have
$$
\sum_{x\in \varphi(\Fp_q)}{\beta(x)}=
\frac{\beta(0)}{2}+\frac{1}{2}\sum_{x\in\Fp_q}{\beta(x^2)}.
$$
\par
If $\beta(x)=\psi(f(x))\chi(g(x))$ where $\psi$ (resp. $\chi$) is an
additive (resp. multiplicative) character, and $f$, $g$ two polynomials,  
the second sum is still of this type, but over the points in $\Fp_q$,
so it can be analyzed in the standard way.
\par
More generally, one can do similar things for any formula of the type
$$
\exists y\ x=h(y)
$$
where $h$ is another polynomial, except that one must handle the
various possibilities for the number of solutions $y$ in
$\Fp_q$ to the equation $h(y)=x$: this kind of step is very clear in
the structure result (Theorem~\ref{th-reduction}).
\end{example}


\begin{example}
We next show that multiplicative characters can be ``tricky'' (as my son
Nicolas would say): consider the formula $\varphi(x)$ in one variable 
$$
\exists y\, x^2+1=y^2
$$
and take $g(x)=x^2+1$ and $\chi$ the multiplicative character of order
$2$ for all finite fields, non-trivial if $p\not=2$. Clearly 
$$
\sum_{x\in\varphi(\Fp_q)}{\chi(g(x))}=\sum_{x\in\varphi(\Fp_q)}{1}
$$
so there is never cancellation. In this case, the reduction in
Theorem~\ref{th-reduction} below applies with $K=1$, $s=0$, $r=0$, $k=1$, 
$\Phi_1(x,y)=\varphi(x)$  and
$h_{1,1}(x,z)=x^2-z^2+1$ (see~(\ref{eq-Phi})), $e=2$, so
$V=V_1=\Aa^1_{\Zz}$, $W=W_{1,1}$ is the affine conic with equation
$X^2-Z^2+1=0$ and $\tau=\tau_{1,1}\,:\, W\ra V$ is given by $X$. The
function $g\circ \tau$ is, on $W$, equal to $X^2+1$, and so is the
square of the function $Z$ on the conic. This is quite obvious, but
indicates that some knowledge of the structure of the definable sets
is required to state a precise criterion for cancellation in a
multiplicative character sum. Still, this knowledge need only be
gained \emph{once} and may then be applied for many different sums. 
\par
Another variant is the following: take the formula in two variables
$$
\varphi(x,y)\,:\, \exists z\ (x^2-y^2+1)^2+z^2=0,
$$
and consider as before the quadratic character sum over
$\varphi(x,y)$ with $g(x)=y^2-1$.
\par
So we introduce the affine variety $W$ with this equation and the
projection $W\ra \Aa^2$. This variety is not absolutely irreducible:
adjoining a square root $\eps$ of $-1$, it splits as
$$
((X^2-Y^2+1)-\eps Z)((X^2-Y^2+1)+\eps Z)=0.
$$
Hence if $-1$ is a square in $\Fp_q$, there is always a value of $z$ for
each $(x,y)$ and the sum is
$$
\sum_{(x,y)\in \Fp_q}{\chi(y^2+1)}=q\sum_{y\in\Fp_q}{\chi(y^2+1)}\ll
q^{3/2},
$$
which gives some cancellation, while on the other hand if $-1$ is not
a square in $\Fp_q$, the points of $W(\Fp_q)$ must have $z=0$, hence
$x^2=y^2-1$ and the sum becomes
$$
\sum_{\stacksum{x\in\Fp_q}{x^2+1\text{ is a square}}}{\chi(x^2)}
$$
which is degenerate as in the previous case, although $g$ is not
the square of a function on $W_{\barre{\Fp}_q}$, but only on the
intersection of its absolutely irreducible components.
\par
Another amusing example is as follows: consider variables $(a,b)$ and
let $\varphi(a,b)$ be the formula
$$
\forall x\ x^2+ax+b\not=0.
$$
\par
For a field $F$ of characteristic $\not=2$, we have $F\models
\varphi(a,b)$ if and only if the discriminant $\Delta(a,b)=a^2-4b$ is
not a square in $F$. Hence if we take $\psi$ trivial, $\chi$ of order
$2$ and $g(a,b)=\Delta(a,b)$, the sum
$$
\sum_{(a,b)\in\varphi(\Fp_q)}{\chi_2(\Delta(a,b))}
$$
has no cancellation. This is only ``trivial'' if one knows about
discriminants of quadratic polynomials, and one can guess that
situations involving invariants of more complicated algebraic forms
will lead to examples which are much less obvious.
\par
Because of this and the lack of current application, we have not tried
to give a cancellation criterion for multiplicative character sum in
this paper, but we hope to investigate this problem further.
\end{example}

\begin{example}
Here is a challenging example: let $n\geq 1$ and consider
the formula $\varphi(a)$, $a=(a_0,a_1,\ldots, a_{n-1})$, which expresses
that the polynomial
\begin{equation}\label{eq-kl-pol}
X^n+a_{n-1}X^{n-1}+\cdots+a_1X+1
\end{equation}
is irreducible and $a_0a_1\cdots a_{n-1}=1$. Consider then the following
exponential sum 
$$
K^*_{n,p}=\sum_{a\in\varphi(\Fp_p)}
{e\Bigl(\frac{a_0+\cdots+a_{n-1}}{p}\Bigr)}
$$
which is thus a subsum of a hyper-Kloosterman $K_{n,p}$ sum in $n$
variables (which only includes the condition $a_0\cdots a_{n-1}=1$ in
the summation). Then we have
\begin{equation}\label{eq-kloos-irred}
|K^*_{n,p}|\ll p^{n-1/2},
\end{equation}
where the implied constant depends on $n$. Compare this with Deligne's
estimate (\cite[Sommes trig. \S 7]{deligne})
$$
|K_{n,p}|\leq (n+1)p^{n/2}.
$$
\par
To prove~(\ref{eq-kloos-irred}), it suffices, according to
Theorem~\ref{th-exp-def}, to show that the function
$$
f(a)=a_0+\cdots +a_{n-1}
$$
is not constant for a positive proportion of
$a\in\varphi(\Fp_p)$. This follows from the following two facts: (1) 
a positive proportion of $a\in \Fp_p^n$ with $a_0\cdots a_{n-1}=1$
define an irreducible polynomial~(\ref{eq-kl-pol}); 
(2) $f$ is not constant on a positive proportion of $a\in \Fp_p^n$
with $a_0\cdots a_{n-1}=1$. Of these, (2) is clear (if one wishes, it
is a consequence of Deligne's estimate for hyper-Kloosterman sums!),
and (1) follows by interpreting the desired number as $1/n$ times the
number of elements in $\Fp_{p^n}^{\times}$ which are of norm $1$ and
are of degree $n$ over $\Fp_p$ (i.e. do not lie in a smaller field).
\par
Here is a natural question: can the estimate for $K^*_{n,p}$ be made
more precise? Is the exponent $n-1/2$ optimal?
\end{example}


\section{Geometric decomposition of exponential sums over definable
sets}
\label{sec-decomp}

Given a formula $\varphi(x,y)$, we start by describing a
``geometric'' decomposition 
of the definable sets $\varphi(\Fp_q,y)$ and any sum over
$\varphi(\Fp_q,y)$, following~\cite{zoe} with some more details. Then
in the next section we will 
apply the sheaf-theoretic and cohomological methods to express the
exponential sums in terms of eigenvalues of the Frobenius operator on
suitable cohomology groups.
\par
It will be noticed that the notation is somewhat involved, since a
large number of parameters occur here and below. In the Appendix, we
give a list of most of them.

\begin{theorem}\label{th-reduction}
Let $\varphi(x,y)$ be a first-order formula in the language of rings
with variables $x=(x_1,\ldots,x_n)$ and parameters $y=(y_1,\ldots,
y_m)$, $n\geq 0$, $m\geq 0$.
\par
There exist the following data, depending only on $\varphi(x,y)$:
\par
\quad \emph{(i)} Integers $K\geq 1$, $s\geq 0$, $e\geq 1$;
\par
\quad \emph{(ii)} A prime power $q_0$;
\par
\quad \emph{(iii)} For $\kappa\leq K$, formulas $\Phi_{\kappa}(x,x',y)$ with
auxiliary parameters $x'=(x'_1,\ldots,x'_{s})$;
\par
\quad \emph{(iv}) For $\kappa\leq K$  and $i\leq e$, affine schemes
$W_{\kappa,i}$ of finite type over $\Zz$ with maps 
$$
\pi_{\kappa,i}\,:\, W_{\kappa,i}\ra \Aa^{s+m}_{\Zz},\quad
\tau_{\kappa,i}\,:\, W_{\kappa,i}\ra \Aa^n_{\Zz};
$$
with the following properties:
\par
\quad \emph{(1)} For every finite field $\Fp_q$ with $q\geq q_0$,
there exists $x'\in \Fp_q^{s}$ such that  
\begin{equation}\label{eq-aux}
\varphi(x,y)\leftrightarrow (\Phi_1(x,x',y)\vee \cdots \vee 
\Phi_K(x,x',y))
\end{equation}
for every $y\in\Fp_q^m$.
\par
\quad \emph{(2)} For every field $F$, every $(x',y)\in F^{s+m}$, the
sets $\Phi_{\kappa}(F,x',y)$, $\kappa\leq K$, are disjoint.
\par
\quad \emph{(3)} For all
finite fields $\Fp_q$ with $q\geq q_0$ and $(x',y)\in\Fp_q^{s+m}$ with
$x'$ chosen so that~\emph{(\ref{eq-aux})} holds, we have 
\begin{equation}\label{eq-imply}
x\in W_{\kappa,i}(\Fp_q)\text{ implies } \tau_{\kappa,i}(x)\in
\Phi_{\kappa}(\Fp_q,x',y)
\end{equation}
and at most $e$ elements of $W_{\kappa,i}(\Fp_q)$ satisfy
$\tau_{\kappa,i}(x)\in\Phi_{\kappa}(\Fp_q,x',y)$.
\par
\quad \emph{(4)} For all finite fields $\Fp_q$ with $q\geq q_0$ and
$(x',y)\in\Fp_q^{s+m}$ with $x'$ chosen so that~\emph{(\ref{eq-aux})} holds, 
and any complex numbers $\beta(x)$ for $x\in\varphi(\Fp_q,y)$, we have
\begin{equation}\label{eq-red-sum}
\sum_{x\in \Phi_{\kappa}(\Fp_q,x',y)}{\beta(x)}=
%
\sum_{1\leq i\leq e}{\frac{(-1)^{i+1}}{i!}
\sum_{\stacksum{x\in W_{\kappa,i}(\Fp_q)}{\pi_{\kappa,i}(x)=(x',y)}}
{\beta(\tau_{\kappa,i}(x))}}.
\end{equation}
\end{theorem}

\begin{proof}
We follow the reduction steps of~\cite[p. 123]{zoe}. This
provides us  with the following data:
\par
\quad $\bullet$ Integers $K\geq 0$, $s\geq 0$, $r\geq 0$, 
$k\geq 0$, $e\geq 1$, and a prime power $q_0$;
\par
\quad $\bullet$ For each integer $\kappa\leq K$, polynomials $f_{\kappa,1}$,
\ldots, $f_{\kappa,r}\in \Zz[X,X',Y]$ where $X'=(X'_1,\ldots, X'_{s})$
is an $s$-tuple of auxiliary parameters;
\par
\quad $\bullet$ For each integer $\kappa\leq K$, polynomials $h_{\kappa,j}
\in \Zz[X,X',Y,Z_j]$ for $1\leq j\leq k$;
\par
\noindent 
which depend only on the formula $\varphi(x,y)$ and have the following
properties: 
\par
\quad $\bullet$ For every finite field $\Fp_q$ with $q\geq q_0$,
the formula $\varphi(x,y)$ is equivalent to the formula
$$
\varphi'(x,y)=\Bigl(\Phi_1(x,x',y)\vee \Phi_2(x,x',y)\vee\cdots\vee
\Phi_K(x,x',y)\Bigr)
$$
where $\Phi_{\kappa}(x,x',y)$ is the formula
\begin{multline}\label{eq-Phi}
\Phi_{\kappa}(x,x',y)=\Bigl(
f_{\kappa,1}(x,x',y)=0\wedge \cdots\wedge f_{\kappa,r}(x,x',y)=0
\\
\wedge (\exists z_1\cdots \exists z_k\ 
h_{\kappa,1}(x,x',y,z_1)=0\wedge \cdots\wedge
h_{\kappa,k}(x,x',y,z_k)=0)
\Bigr),
\end{multline}
where $x'\in \Fp_q^s$ is some value of the auxiliary variables (see
below for a short explanation).
\par
\quad $\bullet$ For any field $F$ and parameters $(x',y)$, the sets
$\Phi_{1}(F,x',y)$,\ldots, $\Phi_K(F,x',y)$ are disjoint in $F^n$;
\par
\quad $\bullet$ For each $\kappa\leq K$ and $(x',y)\in\Fp_q^{s+m}$,
the number of $z=(z_1,\ldots, z_k)\in F^k$ such that $F\models
\Psi_{\kappa}(x,x',y,z)$ is bounded by $e$, where $\Psi_{\kappa}(x,x',y,z)$
is the formula 
\begin{multline} 
\Psi_{\kappa}(x,x',y,z)=\Bigl(
f_{\kappa,1}(x,x',y)=0\wedge \cdots\wedge f_{\kappa,r}(x,x',y)=0
\\
\wedge h_{\kappa,1}(x,x',y,z_1)=0\wedge \cdots\wedge
h_{\kappa,k}(x,x',y,z_k)=0
\Bigr).
\end{multline}
\par
(Precisely, this summarizes the discussion on p. 123 of~\cite{zoe}:
(\ref{eq-Phi}) is the conjunction of the formulas~(1) of loc. cit.,
its precise form given in~(3) in loc. cit.; 
the disjointness of the $\Phi_{\kappa}(\Fp_q,x',y)$ is stated
after~(2) in loc. cit.; the existence and property of $e$ is given
in~(4) in loc. cit.)
\par
Clearly, all this gives already the data of $K$, $s$, $q_0$ and the
auxiliary formulas $\Phi_{\kappa}(x,x',y)$ such that
(1) and (2) of the theorem hold. Before continuing to the second part,
we explain the  occurrence of the auxiliary parameters $x'$: they
correspond to extra symbols $c_i$ in the language of enriched fields
discussed in~\cite[\S 2]{zoe}. For any field $F$, those constants are
interpreted as coefficients of some irreducible monic polynomial in
$F[T]$. To see the relevance with definable sets $\varphi(\Fp_q,y)$,
notice that for instance the universal formula
\begin{equation}\label{eq-universal} 
\forall t\ a_0+a_1t+a_2t^2\not=0
\end{equation}
can be restated, for a finite field $\Fp_q$, by the existential
formula stating that $a_0+a_1t+a_2t^2$ splits in linear factors over
the field $\Fp_{q^2}$, with no root in $\Fp_q$. This can be expressed
in terms of the coefficients 
of an irreducible monic polynomial $f_2=c_0+c_1T+T^2$ by factorizing
$a_2(T-y_1)(T-y_2)$ and then expressing $\Fp_q$-rationally the
equality  
$$
a_0+a_1t+a_2t^2=a_2(t-y_1)(t-y_2)\text{ with } y_1,\ y_2\in 
\Fp_{q^2}-\Fp_q
$$
in the basis
$(1,\alpha_2)$, where $\alpha_2$ is a root of $f_2$.\footnote{\
Alternately,~(\ref{eq-universal}) can be approached by Galois theory;
this would lead to the use of the Galois stratification method
described in~\cite{fried-haran-jarden}.} It is also useful to mention
that it is in constructing $\Phi_{\kappa}$ that the restriction to
large enough finite fields is introduced, this coming from a result of
van den Dries, based on the Lang-Weil estimate and some model theory
(see~\cite[2.4]{zoe}). 
\par
Now, to simplify notation for the proof of (3) and (4) of the
theorem, which concern the formulas $\Phi_{\kappa}(x,x',y)$
individually, we drop the subscript $\kappa$ and we incorporate the 
new parameters $x'$ into $y$ (so that $s+m$ is now denoted $m$). 
\par
Let $V$ denote the zero set of the polynomials $f_{1}$, \ldots,
$f_r\in \Zz[X,Y]$, seen as a closed subscheme of
$\Aa_{\Zz}^{n+m}$.  We have the projection $V\ra \Aa_{\Zz}^m$
given by the parameter $y$, and we denote by $V_y$ the fibers.
Let $W$ denote the common zero set in $\Aa_{\Zz}^{n+m+k}$ of the
polynomials $f_i$ and the polynomials $h_1$, \ldots, $h_k$.
Denote by $\pi\,:\, W\ra V$ the obvious projection.
\par
For $j\geq 1$, denote
$$
W_y(\Fp_q)_j=\{x\in \Phi(\Fp_q,y)\,\mid\,
|\pi^{-1}(x,y)\cap W(\Fp_q)|=j\}.
$$
\par
Let now $\beta(x)$ be arbitrary complex numbers
defined for $x\in \Phi(\Fp_q,y)$. 
By~(\ref{eq-Phi}), for any $x\in \Phi(\Fp_q,y)$, we have $(x,y)\in
V(\Fp_q)$ and $\pi^{-1}(x,y)\cap W(\Fp_q)\not=\emptyset$, so that
$$
\sum_{x\in \Phi(\Fp_q,y)}{\beta(x)}
=\sum_{j\geq 1}{\sum_{x\in W_y(\Fp_q)_j}{
\beta(x)}}
$$
and from the defining property of $e$, this reduces to
\begin{equation}\label{eq-red}
\sum_{x\in \Phi(\Fp_q,y)}{\beta(x)}=\sum_{1\leq j\leq e}{
\sum_{x\in W_y(\Fp_q)_j}{
\beta(x)}}.
\end{equation}
\par
To deal with the fact that the fibers of $\pi$ do not
necessarily all have the same cardinality, we now use the same 
combinatorial procedure as in~\cite[p. 124]{zoe}, although we make the
result more explicit.
\par
For  $1\leq j\leq e$, let
$W_j$ denote the intersection in $\Aa_{\Zz}^{n+m+jk}$ of the $j$-fold
fiber product of $W$ over $V$ with the open
subscheme  $U_j$ consisting of points $(x,y,z_1,\ldots, z_j)$ where
all the $z_i$ are distinct $k$-tuples, i.e. in point terms we have
$$
W_j(A)=U_j(A)\cap \{(x,y,z_1,\ldots, z_j)\,\mid\, (x,y,z_i)\in W(A)\text{
  for } 1\leq i\leq k\}
$$
for any ring $A$ (this corresponds to the formula
$\Psi_j(X,Y,Z^1,\ldots, Z^j)$ of loc. cit.) On $W_j$ we have the
maps
$$
\pi_j\ \begin{cases}
\quad W_j \ra \Aa_{\Zz}^m\\
(x,y,z_1,\ldots, z_j)\mapsto y
\end{cases}
\quad\quad\text{ and }
\quad\quad
\tau_j\ 
\begin{cases}
\quad W_j\ra \Aa^n_{\Zz}\\
(x,y,z_1\ldots,z_j)\mapsto x;
\end{cases}
$$
we will denote by $W_{j,y}$ the fiber of $\pi_j$ over $y$. All this
(with the omitted dependency on $\kappa$) gives the data (iv) of the
Theorem and (3) is satisfied by construction.
\par
Next we claim that the following combinatorial formulas hold: denoting
$$
(i)_j=i(i-1)\cdots (i-j+1)=j!\binom{i}{j},
$$
we have for $1\leq j\leq e$
\begin{equation}\label{eq-incl-excl}
\sum_{i=j}^e{(i)_j\sum_{x\in W_y(\Fp_q)_i}{\beta(x)}}=
\sum_{(x,y,z_1,\ldots, z_j)\in W_j(\Fp_q)}{\beta(x)}=
\sum_{x\in W_{j,y}(\Fp_q)}{\beta(\tau_j(x))}.
\end{equation}
Indeed, for each $x\in W_y(\Fp_q)_i$, the set $\pi^{-1}(x,y)\cap
W(\Fp_q)$ has $i$ elements and for $j\leq i\leq e$, any
ordered subset of length $j$ naturally gives a point of
$W_{j,y}(\Fp_q)$. All of these points are distinct, and of course
there are precisely $(i)_j$ of them. 
\par
Now we use the following elementary lemma (presumably standard in
combinatorics): 

\begin{lemma}
Let $e\geq 1$. Let $x_j$, $y_i$ be complex numbers defined for $1\leq
i,j\leq e$, such that 
\begin{equation}
\label{eq-triang-sys}
\sum_{i=j}^e{(i)_j x_{i}}=y_j,
\end{equation}
for $1\leq j\leq e$. Then we have
$$
\sum_{1\leq i\leq e}{x_i}=
\sum_{1\leq j\leq e}{\frac{(-1)^{j+1}}{j!}y_j}.
$$
\end{lemma}

\begin{proof}
To see this, solve first the triangular system~(\ref{eq-triang-sys})
using $(i)_j=i!/(i-j)!$ to get
$$
x_i=\sum_{j=i}^{e}{\frac{(-1)^{i+j}}{i!(j-i)!}y_j},
$$
(as easily checked using the binomial theorem), and then sum over $i$
to get 
\begin{align*}
\sum_{1\leq i\leq e}{x_i}&=\sum_{1\leq i\leq e}{
\sum_{i\leq j\leq e}{\frac{(-1)^{i+j}}{i!(j-i)!}y_j}}
=
\sum_{1\leq j\leq e}{(-1)^jy_j
\Bigl(\sum_{1\leq i\leq j}{\frac{(-1)^i}{i!(j-i)!}}\Bigr)}\\
&=\sum_{1\leq j\leq e}{\frac{(-1)^{j+1}}{j!}y_j}
\end{align*}
by the binomial theorem again.
\end{proof}

Applied to 
$$
x_i=\sum_{x\in W_y(\Fp_q)_i}{\beta(x)},\quad 
y_j=\sum_{x\in W_{j,y}(\Fp_q)}{\beta(\tau_j(x))},
$$
we get by~(\ref{eq-red}) and~(\ref{eq-incl-excl}) that
$$
\sum_{x\in\Phi(\Fp_q,y)}{
\beta(x)}=
\sum_{1\leq j\leq e}{\frac{(-1)^{j+1}}{j!}
\Bigl(\sum_{x\in W_{j,y}(\Fp_q)}{\beta(\tau_j(x))}\Bigr)},
$$
which is the final conclusion~(\ref{eq-red-sum}), taking account the
change of notation made before.
\end{proof}

\section{Estimates arising from the decomposition of definable sets}
\label{sec-estimates}

For the second part of the reduction of exponential sums over
definable sets, recall that given a prime power $q$ and an integer
$k\geq 0$, a complex number $\alpha\in\Cc$ is a $q$-Weil number of weight  
$k$ if $\alpha$ is an algebraic integer, and any 
Galois-conjugate $\alpha'$ of $\alpha$ (i.e., any root of the minimal
polynomial $P\in \Zz[T]$ of $\alpha$) satisfies $|\alpha'|=q^{k/2}$.
It will be convenient to call a \emph{signed} $q$-Weil number a pair $(\pm
1,\alpha)$ of a sign and a $q$-Weil number. Usually we just write
$\alpha$ and the sign is denoted $\eps(\alpha)$. When a Weil number is
written down explicitly as a complex number and claimed to be a signed
Weil number, this means that the sign is $+1$. For instance, this
applies to $\alpha=q^w$, a $q$-Weil number of weight $2w$.
\par
To compute and estimate the exponential sums $S(y,\varphi,\Fp_q)$, we
will apply the Grothendieck-Lefschetz trace formula and the Riemann
Hypothesis over finite fields proved by Deligne to the auxiliary
varieties $\pi_{\kappa,i}^{-1}(x',y)\subset W_{\kappa,i}$ arising from
Theorem~\ref{th-reduction}. This naturally falls in two steps: first,
we examine the number of points of summation, then we analyze when an
exponential sum exhibits cancellation.
\par
We will first perform both steps for a single $W_{\kappa,i,y}$, i.e.,
for a ``classical'' exponential sum (in logical terms, one
corresponding to a positive quantifier free formula). The only
subtlety is that we do 
not know if $W_{\kappa,i,y}$ is absolutely irreducible or not (the
case usually treated in the literature), which affects the precise
counting of points and the non-triviality of the estimates.

\begin{proposition}\label{pr-variety}
Let $W\subset \Aa^{n+m}_{\Zz}$ be an affine subscheme, let $q$ be a power
of the prime $p$, $(\psi,\chi,f,g)$ data defining the exponential sum 
$$
S(y,W,\Fp_q)=\sum_{(x,y)\in W(\Fp_q)}{\psi(f(x))\chi(g(x))}
$$
for $y\in \Fp_q^m$ over $W$. 
\par
\emph{(i)} There exist an integer $B_1\geq 0$,
depending only on $W$, and for all $y\in\Fp_q^m$, there exists an
integer $\delta(y)$, $0\leq \delta(y)\leq 
n$, and signed $q$-Weil numbers $\alpha_j(y)$ for $1\leq j\leq
\beta\leq B_1$, where $\beta$ may depend on $y$, such that 
$$
w(\alpha_{j}(y))\leq 2n,\quad
\max w(\alpha_{j}(y))=2\delta(y),\text{ and }
\alpha_{j}(y)=q^{\delta(y)}\text{ if }
w(\alpha_{j}(y))=2\delta(y),
$$
and
\begin{equation}\label{eq-number-v}
\sum_{(x,y)\in W(\Fp_q)}{1}=\sum_{1\leq j\leq \beta}{
\eps(\alpha_j(y))\alpha_j(y)}.
\end{equation}
\par
\emph{(ii)} There exist an integer $B_2\geq 0$,
depending only on $W$ and the degree of $f$ and $g$, and for all
$y\in\Fp_q^m$ there exist $q$-Weil
numbers $\beta_j(y)$, for $1\leq j\leq \gamma\leq B_2$, where 
$\gamma$ may depend on $y$, such that
$$
w(\beta_{j}(y))\leq 2\delta(y)\text{ for all $j$},\\
$$
and 
\begin{equation}\label{eq-sum-v}
\sum_{(x,y)\in W(\Fp_q)}{\psi(f(x))\chi(g(x))}=
\sum_{1\leq j\leq \gamma}{\eps(\beta_j(y))\beta_j(y)}.
\end{equation}
Moreover, there exists $\eta>0$, depending only on $W/\Zz$ such that
if $p$ is large enough, $p\geq p_0$ where $p_0$ depends only  
on $W$ and the degrees of $f_1$ and $f_2$, we have
$w(\beta_j(y))<2\delta(y)$ for all $j$ unless either 
$\psi$ is trivial or there exists $c\in\Fp_q$ such that 
\begin{equation}\label{eq-pos-prop}
\sum_{\stacksum{(x,y)\in W(\Fp_q)}{f(x)=c}}{1}\geq \eta
\sum_{(x,y)\in W(\Fp_q)}{1}.
\end{equation}
\end{proposition}

\begin{proof}
We denote by $W_y$ the fiber over $y$ of $W/\Fp_q$. For reasons that
will become clear soon, we first replace $W_y$ by 
the subscheme of $\Aa^{n+m}_{\Fp_q}$ obtained by performing the
decomposition-intersection process (described in~\cite[\S 1]{zoe}),
i.e., we replace $W_y$ by $V_y$ such that
\par
\quad (i) $W_y(\Fp_q)=V_y(\Fp_q)$;
\par
\quad (ii) The absolutely irreducible components of $V_y$ are defined
over $\Fp_q$.
\par
It follows from~\cite[Pr. 1.7]{zoe} that $V_y$ can be defined as the
zero set of $N$ polynomials in $\Fp_q[X_1,\ldots, X_n]$ of degree
$\leq E$, and moreover $V$
has at most $I$ absolutely irreducible components, where $N$, $E$
and $I$ depend only on $W/\Zz$ (in particular are independent of $y$).
\par
By~(i), we have for all $y\in\Fp_q^m$
$$
S(y,W,\Fp_q)=S(V_y,\Fp_q).
$$
Note it is possible that $V_y=\emptyset$ (for instance for $X^2+1=0$
if $-1$ is not a square in $\Fp_q$); if that is the case, we can take
$B_1=B_2=0$, and we put $\delta(y)=0$. So assume $V_y$ is not empty.
\par
Fix a prime $\ell\not=p$, where $p$ is the characteristic of
$\Fp_q$. The formalism of the Lang torsor (see
e.g.~\cite[Sommes trig. 1.4]{deligne},~\cite[4.3]{gsksm} 
or the sketch in~\cite[11.11]{ant}) provides us with the
lisse $\barre{\Qq}_{\ell}$-adic sheaf of rank $1$
$$
\sheaf{L}=\sheaf{L}_{\psi(f)}\otimes \sheaf{L}_{\chi(g)}
$$ 
on $V_y$ (which depends on
$p$, $\psi$, $\chi$, $f$ and $g$) such that  the
local trace of a geometric Frobenius element $\frob_{x,\Fp_q}$
at a rational point $x\in V_y(\Fp_{q})$ is given by 
\begin{equation}\label{eq-local-trace}
\Tr(\frob_{x,\Fp_{q}}\,\mid\, \mathcal{L})=
\psi(f(x))\chi(g(x))
\end{equation}
and therefore
$$
S(V_y,\Fp_q)=
\sum_{x\in V_y(\Fp_{q})}{
\Tr(\frob_{x,\Fp_{q}}\,\mid\, \mathcal{L})}.
$$
\par
Let $\barre{V}_{y}$ denote the base
change of $V_y$ to the algebraic closure of $\Fp_q$. The
Grothendieck-Lefschetz trace formula  (see e.g.~\cite[Rapport]{deligne})
gives the cohomological expression
$$
S(V_y,\Fp_q)=\sum_{i\geq 0}{(-1)^i \Tr(\frob\,\mid\, 
H_c^i(\barre{V}_y,\sheaf{L}))}
$$
where $H_c^i$ are the compactly supported $\ell$-adic cohomology
groups with coefficient in the sheaf $\sheaf{L}$ and $\frob$ denotes
the geometric Frobenius operator over $\Fp_q$ which acts naturally on
those finite dimensional $\barre{\Qq}_{\ell}$-vector spaces. Of course, the
sum is in fact finite, and more precisely the space $H^i_c$ vanishes
for $i>2\delta(y)$, where $\delta(y)$ is the maximal dimension of
an irreducible component of $\barre{V}_y$. Because $V_y$ is a subset
of affine $n$ space, we have $\delta(y)\leq n$.
\par
By Deligne's fundamental result (his far-reaching generalization of
the Riemann Hypothesis), since the sheaf $\sheaf{L}$ is punctually
pure of weight $0$ (its local traces being roots of unity), the
eigenvalues of $\frob$ acting on
$H_c^i(\barre{V}_y,\mathcal{L})$ are $q$-Weil numbers of
weight $\leq i$ (see~\cite[Th. 3.3.1]{weilii}).
\par
Thus we obtain~(\ref{eq-number-v}) in the case $\mathcal{L}=\Qq_{\ell}$
with the family of $q$-Weil numbers $\alpha_{j}(y)$ being the
eigenvalues (indexed by $j$, and with multiplicity) of the geometric
Frobenius on all the 
cohomology groups  $H^i_c(\barre{V}_y,\Qq_{\ell})$ for
$i\leq 2\delta(y)$, those becoming ``signed'' by the factor
$(-1)^i$. 
\par
In this case still, for $i=2\delta(y)$, let $\barre{U}_y$ be an
irreducible component of $\barre{V}_y$; by (ii) above,
$\barre{U}_y=U_y\times\barre{\Fp}_q$ for some irreducible
$U/\Fp_q$. By Poincaré duality (see~\cite[Sommes trig.,
Rem. 1.18d]{deligne}), we have 
$$
H^{2\delta(y)}_c(\barre{U}_y,\Qq_{\ell})=
(\Qq_{\ell})_{\pi_1(\barre{U}'_y,\barre{\eta})}(-2\delta(y))
=(\Qq_{\ell})(-2\delta(y))
$$
where $U'_y\subset U_y$ is a smooth, dense open subscheme of $U_y$,
and $\pi_1(\barre{U}'_y,\barre{\eta})$ is the geometric fundamental
group of $\barre{V}_y$ with respect to some geometric point
$\barre{\eta}$. This precisely means that $H^{2\delta(y)}_c$ is of
dimension $1$ and the geometric Frobenius $\frob$ of $\Fp_q$ acts by
multiplication by $q^{\delta(y)}$. (The only delicate point is that if
we use an irreducible component not defined over $\Fp_q$, it is the
Frobenius $\frob^{\nu}$ of a field on which it is defined that acts by
multiplication  by $q^{\nu\delta(y)}$).
\par
Since moreover it is known that
$H^{2\delta(y)}_c(\barre{V}_y,\Qq_{\ell})$ is the direct sum of the
corresponding cohomology of the irreducible components, 
it follows that all eigenvalues of weight $2\delta(y)$  are equal to
$q^{2\delta(y)}$ and that the multiplicity is the number of irreducible
components of dimension $2\delta(y)$ (see~\cite[Sommes trig.,
Remarques 1.18 (d)]{deligne}. This gives all the properties of 
$\alpha_j$ stated in the first part of Proposition~\ref{pr-variety}.
\par
Similarly in the general case we obtain~(\ref{eq-sum-v}) with the
analogue eigenvalues for $H^i_c(\barre{V}_y,\mathcal{L})$,
and the weight of these is $\leq 2\delta(y)$ as stated.
\par
A crucial point that remains to be checked is that the total number of
eigenvalues (i.e. the numbers denoted $\beta$ and $\gamma$ in the
statement of the Proposition) is indeed bounded by $B_1$ or $B_2$
depending only on $W$ and (in the case of $B_2$) on
the degrees of $f_1$, $f_2$, $g_1$, $g_2$.
\par
This is precisely given by a very useful result of
Katz~\cite[Th. 12]{katz}, because we can 
embed $\barre{V}_y$ as a closed subscheme of an affine space
$\Aa^M_{\barre{\Fp}_q}$ where $M$ depends only on $W/\Zz$, using equations
of degree and number bounded only in terms of $W/\Zz$ (by the
uniformity of the decomposition-intersection procedure described at
the beginning).\footnote{\ It may seem surprising, but it is indeed
  true that $B_2$ is independent of the order of the multiplicative
  character $\chi$.}
\par
Now we analyze when we can get some cancellation in~(\ref{eq-sum-v}). 
We assume that $\psi$ is non-trivial and that
$\max w(\beta_j(y))=2\delta(y)$ and will show then that~(\ref{eq-pos-prop})
holds for some $c$.
\par
The hypothesis implies that for some irreducible component
$\barre{U}_y=U_y\times\barre{\Fp}_q$ of $\barre{V}_y$, the cohomology
group $H^{2\delta(y)}_c(\barre{U}_y,\mathcal{L})$ does not vanish
(since all $H^i_c$ with $i>2\delta(y)$ do vanish and those with
$i<2\delta(y)$ yield Weil numbers with smaller weight).
\par
If need be, we replace $U_y$ without changing notation by a smooth
non-empty open subscheme. Since $\mathcal{L}$ restricted to
$U_y$ is a lisse sheaf on a smooth connected scheme, we
have by Poincaré duality (as before) the co-invariant formula
\begin{equation}\label{eq-coinv}
H^{2\delta(y)}_c(\barre{U}_y,\mathcal{L})=
(\mathcal{L}_{\eta})_{\pi_1(\barre{U}_y,\eta)}(-2\delta(y))
\end{equation}
for any geometric point $\eta$ of $\barre{U}_y$, $\mathcal{L}_{\eta}$
being the fiber of $\mathcal{L}$ at $\eta$. Since $\mathcal{L}$ is of
rank $1$ and the cohomology does not vanish, we must therefore have
$$
(\mathcal{L}_{\eta})_{\pi_1(\barre{U}_y,\eta)}(-2\delta(y))
=\mathcal{L}_{\eta}(-2\delta(y)), 
$$
i.e., the sheaf $\mathcal{L}$ is geometrically trivial.
From the exact sequence
$$
1\ra \pi_1(\barre{U}_y,\eta)\ra \pi_1(U_y,\eta)\ra 
\Gal(\barre{\Fp}_q/\Fp_{q})\ra 1,
$$
and~(\ref{eq-coinv}), it follows that $\mathcal{L}$ on $U_y$ ``is'' a
character of $\pi_1(U_y,\eta)$ which comes from a character of
$\Gal(\barre{\Fp}_q/\Fp_{q})$. This means in particular that
the local trace of $\mathcal{L}$ at the geometric Frobenius of a point
of $U_y(\Fp_{q})$ is constant 
i.e. (by~(\ref{eq-local-trace})), for all $x,x_0\in
U_y(\Fp_{q})$ we have
$$
\psi(f(x))\chi(g(x))=
\psi(f(x_0))\chi(g(x_0)).
$$
Taking the $D$-th power (where $D$ is the order of $\chi$) yields
the equality 
\begin{equation}\label{eq-constant-coset}
\psi(D(f(x)-f(x_0)))=1.
\end{equation}
\par
We can write $\psi(x)=e(\Tr(ax)/p)$ for some $a\in\Fp_q^{\times}$,
since $\psi$ is assumed non-trivial. Then~(\ref{eq-constant-coset})
and the injectivity of $x\mapsto e(x/p)$ on $\Fp_p$ mean that for
all $x\in U_y(\Fp_{q})$, the element
$aDf(x)$ is in a coset of the kernel of the trace from
$\Fp_{q}$ to $\Fp_p$. Such a coset has cardinality
$p^{-1}|\Fp_{q}|$, but on the other hand 
if $f$ is non-constant on $U_y$, and if $q$ is large
enough, there are at least $q/\max(\deg(f_1),\deg(f_2))$ elements of the
form $aDf(x)$ in $\Fp_{q}$. (Note $D\equiv 1\mods{p}$.)
\par
So under the hypothesis that $w(y)=2\delta(y)$,
we see that if $p$ is large enough (depending on $W$ and the degrees
of $f_1$ and $f_2$) it must be the case that
$f|\barre{U}_y$ is constant. Since the number $\mu(y)$ of 
irreducible components of $V_y$ of dimension $\delta(y)$ is bounded in
terms of $W/\Zz$ only, the standard counting 
$$
|W_y(\Fp_q)|=\mu(y)q^{\delta(y)}+O(q^{\delta(y)-1/2})
$$
(with absolute constant depending only on $W/\Zz$) coming
from~(\ref{eq-number-v}) and its analogue
$$
|U_y(\Fp_q)|=q^{\delta(y)}+O(q^{\delta(y)-1/2})
$$
show that~(\ref{eq-pos-prop}) holds for $c$ the constant value of $f$
on $U_y(\Fp_q)$, $p$ large enough in terms of $W/\Zz$.
(We could also have argued more geometrically using the
triviality of the Artin-Schreier covering associated to
$\mathcal{L}$, as in~\cite[p. 24, last \S]{deligne}).
\end{proof}

\begin{remark}
There is presumably a cohomological interpretation of the
decomposition-inter\-section process, which means essentially a general
analysis of the cohomology of $W/\Fp_q$ with coefficient (at least) in
a sheaf of the type $\sheaf{L}$ above, in the case where $W$ is not
absolutely irreducible, explaining how the trace formula boils down to
that of the variety obtained by the process. However the author has
not found such a result in the literature. 
\par
A drawback of using this procedure is that the Weil numbers involved
are not uniquely determined by the exponential sum $S(y,W,\Fp_q)$ and
its companions
\begin{equation}\label{eq-companions}
S_{\nu}(y,W,\Fp_{q})=
\sum_{(x,y)\in W(\Fp_{q^{\nu}})}{\psi(\Tr f(x))\chi(Ng(x))}
\end{equation}
(where $\Tr$ and $N$ are the norms from $\Fp_{q^{\nu}}$ or
$\Fp_{q^{\nu}}^{\times}$ to $\Fp_{q}$), 
rather they are uniquely determined by the sequence of sums
$$
S_{\nu}(V_y,\Fp_{q^{\nu}})=\sum_{x\in V_y(\Fp_{q^{\nu}})}{
\psi(\Tr f(x))\chi(Ng(x))}=
\sum_{1\leq j\leq \gamma}{\eps(\beta_j(y))
\beta_j(y)^{\nu}}
$$
for $\nu\geq 1$, 
in the sense that the
multiplicity (with signs taken into account)
$$
\sum_{\beta_j(y)=\alpha}{\eps(\beta_j(y))}
$$
of any $q$-Weil number $\alpha$ is determined by
$(S_{\nu}(V_y,\Fp_q))$. However, those sums for $\nu\geq 2$ are not
necessarily related to the original sum $S(y,W,\Fp_q)$.
\par
One can ignore all this, in a sense, and apply the
Grothendieck-Lefschetz trace formula and Deligne's Theorem to $W_y$
directly. But then the computation of the eigenvalues for the
top-dimensional cohomology is not 
valid unless $\dim W_y/\Fp_q=\delta(y)$; note that this will often be
the case for a concrete $W/\Zz$ (e.g. if $W/\Zz$ is generically
absolutely irreducible), but in our applications to definable sets,
the auxiliary varieties are not so well controlled, especially as they
depend on the value of the auxiliary parameters $x'$.
\end{remark}

We come back to a general formula $\varphi(x,y)$ and start with the
special case $f=0$, $g=1$, $\psi=1$, $\chi=1$ that 
counts the points of the definable sets.

\begin{theorem}\label{th-weil-numbers-count}
Let $\varphi(x,y)$ be a first-order formula in the language of rings
and let $K$ and $e$ be given by Theorem~\ref{th-reduction} for
$\varphi$. There exists an integer $B_1\geq 0$, depending
only on $\varphi$, with the following property: 
for all $q$ large enough and all $y\in\Fp_q^m$,
there exist signed $q$-Weil numbers
$$
\alpha_{\kappa,i,1}(y),\ldots,\alpha_{\kappa,i,\beta}(y),\text{ with }
1\leq \kappa\leq K,\quad 1\leq i\leq e,
$$
with $\beta\leq B_1$, $\beta$ possibly depending on $y$, such that
\begin{gather*}
w(\alpha_{\kappa,i,j}(y))\leq 2n,\quad 
\max w(\alpha_{\kappa,i,j}(y))=2\delta(y)\text{ is even},\\
\alpha_{\kappa,i,j}(y)=q^{\delta(y)}\text{ if }
w(\alpha_{\kappa,i,j}(y))=2\delta(y),
\end{gather*}
and we have
\begin{equation}\label{eq-number}
|\varphi(\Fp_{q},y)|=\sum_{1\leq \kappa\leq K}{\sum_{1\leq i\leq e}{
\frac{(-1)^{i+1}}{i!}
\sum_{1\leq j\leq \beta}{\eps(\alpha_{\kappa,i,j}(y))
\alpha_{\kappa,i,j}(y)}}}.
\end{equation}
\par
Finally, if $\varphi(\Fp_q,y)$ is not empty, the multiplicity ``up to
sign'' of  $\alpha=q^{2\delta(y)}$ is strictly positive, i.e.
$$
\mu(y)=\triplesum_{\stacksum{\kappa,\ \ i,\ \ j}
{\alpha_{\kappa,i,j}(y)=q^{2\delta(y)}}}{\frac{(-1)^{i+1}
\eps(\alpha_{\kappa,i,j}(y))}{i!}}>0
$$
\end{theorem}

We thus recover the main theorem of~\cite{zoe}.

\begin{corollary}[Chatzidakis-van den Dries-Macintyre]\label{cor-number}
\emph{(1)} For all finite fields $\Fp_q$ and $y\in \Fp_q^m$, we have
$$
|\varphi(\Fp_q,y)|=\mu(y) q^{\delta(y)}
+O(q^{\delta(y)-1/2})
$$
where
$$
\mu(y)=\triplesum_{\stacksum{\kappa,i,j}
{w(\alpha_{\kappa,i,j}(y))=2\delta(y)}}
{\frac{(-1)^{i+1}}{i!}}\in \Qq
$$
is $>0$ unless $\varphi(\Fp_q,y)=\emptyset$, and the implied constant
depends only on $\varphi$. 
\par
\emph{(2)} There exist only finitely many pairs $(d,\mu)$ with $d\geq
0$ and $\mu\in\Qq$ which arise as $(\delta(y),\mu(y))$ for some finite
field $\Fp_q$ and $y\in\Fp_q^m$. 
\par
\emph{(3)} For each pair $(d,\mu)$, $d\geq 0$ an integer and
$\mu\in\Qq$, that can arise as $(\delta(y),\mu(y))$, there exists a
formula $\mathcal{C}_{d,\mu,\varphi}$ in the language of rings,
depending only on $d$, $\mu$ and $\varphi$, such that for any finite
field $\Fp_q$, we have $\Fp_q\models \mathcal{C}_{d,\mu,\varphi}$ if
and only if $(\delta(y),\mu(y))=(d,\mu)$.
\end{corollary}

\begin{proof}
(1) For $q$ large enough, this is obvious from~(\ref{eq-number}) (with
$\nu=1$) and the stated properties of the weights of the
$\alpha_{\kappa,i,j}(y)$, recalling that the sign
$\eps(q^{2\delta(y)})$ is $+1$ by convention; in fact we get
$$
\Bigl||\varphi(\Fp_q,y)|-\mu(y) q^{\delta(y)}\Bigr|
\leq 3KB_1q^{\delta(y)-1/2}
$$
for $q\geq q_0$ large enough so that~(\ref{eq-number})
applies. 
\par
Replacing the constant $3KB_1$ by $\max(3KB_1,C)$ where $C\geq 0$
satisfies 
$$
\Bigl||\varphi(\Fp_q,y)|-\mu(y) q^{\delta(y)}\Bigr|\leq Cq^{\delta(y)-1/2}
$$
for $q\leq q_0$, we obtain the result for all $q$. 
\par
(2) This is clear since we have the trivial bound
$|\varphi(\Fp_q,y)|\leq q^n$ for all $n$ which gives $\delta(y)\leq
n$, and because there are at most $2^{K+e+B_1}$ choices of subsets of
summation of indices $(\kappa,i,j)$ that can occur in defining
$\mu(y)$. 
\par
(3) This is proved in~\cite[Prop. 3.8]{zoe} (and can be guessed from
the proof of Theorem~\ref{th-reduction} and
Proposition~\ref{pr-variety}).
\end{proof}

We refer to~\cite[4.9,4.10]{zoe} for intrinsic interpretations of the
``dimension'' $\delta(y)$ and the ``measure'' $\mu(y)$ in the context
of pseudo-finite fields.
\par
We can now come back to the general exponential
sums~(\ref{eq-exp-def}) as we are in a position to compare the
estimates obtained with the number of points of summation.

\begin{theorem}\label{th-weil-numbers-exps}
Let $(f_1,g_1,f_2,g_2,\{\psi\},\{\chi\})$ be the data
defining a family of exponential sums~\emph{(\ref{eq-exp-def})} over
the definable sets $\varphi(\Fp_q,y)$. 
\par
\emph{(1)}
There exists an integer $B_2\geq 0$, depending
only on $\varphi$ and the degrees of $f_1$, $g_1$, $f_2$, $g_2$, with
the following property: 
for all $q$ large enough and all $y\in\Fp_q^m$,
there exist signed $q$-Weil numbers
$$
\alpha_{\kappa,i,1}(y),\ldots,\alpha_{\kappa,i,\beta}(y),\text{ with }
1\leq \kappa\leq K,\quad 1\leq i\leq e,
$$
for $\beta\leq B_2$, $\beta$ depending possibly on $y$, such that
\begin{gather*}
w(\alpha_{\kappa,i,j}(y))\leq 2\delta(y),\text{ for all $\kappa$, $i$, $j$}
\end{gather*}
and  
\begin{equation}\label{eq-decomp}
S(y,\varphi,\Fp_q)=\sum_{1\leq \kappa\leq K}{\sum_{1\leq i\leq e}{
\frac{(-1)^{i+1}}{i!}
\sum_{1\leq j\leq \beta}{\eps(\alpha_{\kappa,i,j}(y))
\alpha_{\kappa,i,j}(y)}}},
\end{equation}
hence
$$
|S(y,\varphi,\Fp_q)|\leq 3KB_2q^{w(y)/2}=3K_2Bq^{\delta(y)-\gamma(y)/2}.
$$
\par
\emph{(2)} Let 
$$
w(y)=\max_{\kappa,i,j} w(\alpha_{\kappa,i,j}(y))\leq 2\delta(y),
$$
denote the maximal weight of the Weil numbers occurring in this
decomposition, and $\gamma(y)=2\delta(y)-w(y)\geq 0$. 
\par
There exists $\eta>0$ depending only on $\varphi(x,y)$ such that for
$p$ is large enough, depending only on $\varphi(x,y)$ and the degrees
of $f_1$ and $f_2$, we have
$w(y)<2\delta(y)$, i.e., $\gamma(y)>0$, unless $\psi$ is trivial or
there exists some $c\in\Fp_q$ with
$$
\sum_{\stacksum{x\in \varphi(\Fp_q,y)}{f(x)=c}}{1}\geq \eta
\sum_{x\in \varphi(\Fp_q,y)}{1}.
$$
\end{theorem}

\begin{remark}
(1) We emphasize that in part (2), we must have $p$ large enough, and not
only $q$. This is necessary even for classical sums, since if we fix
the additive characters by defining $\psi$ on $\Fp_{p^{\nu}}$ as
$e(\Tr(x)/p)$, we have $\psi(x^p-x)=1$ for all $x\in \Fp_{p^{\nu}}$,
so that for a polynomial $f$ congruent modulo $p$ to a polynomial of
the form $g(x)^p-g(x)$, we have
$$
\sum_{x\in\Fp_{p^{\nu}}}{\psi(f(x))}=p^{\nu}
$$
for all $\nu\geq 1$. And of course, such  congruences can a priori
hold for a large number of distinct primes. So our statement only
provides a criterion for cancellation in the ``horizontal'' direction
$p\ra +\infty$ which is of most interest in analytic number theory. 
\par
(2) Note that of course it is only when $\gamma(y)>0$ that this has any
interest. The condition stated to ensure this does not sound
particularly convenient, but that is partly because of the generality
allowed. Essentially, for additive character sums, it says that there
is some cancellation in 
$S(y,\varphi,\Fp_q)$, unless it turns out that $f$ is constant for a
positive proportion of the points of 
summation. This is a generic property;  one may say
that it corresponds to $\varphi(\Fp_q,y)$ being ``transverse'' in some
sense to the level sets of $f$. In concrete applications, it should be
clearer how to check this condition. 
\end{remark}

\begin{proof}[Proof of Theorem~\ref{th-weil-numbers-count}
and Theorem~\ref{th-weil-numbers-exps}]
By~(\ref{eq-red-sum}) with $\beta(x)=\psi(f(x))\chi(g(x))$, we are
reduced to the sums over the 
$\Fp_q$-rational points of the fibers of $\pi_{\kappa,i}\,:\,
W_{\kappa,i}\ra \Aa^n$, which are algebraic varieties. We consider
each in turn, assuming $q\geq q_0$ as in Theorem~\ref{th-reduction} 
and always consider that a choice of the auxiliary parameters $x'$ has
been performed, which we incorporate into $y$ for clarity. In any case
we apply Proposition~\ref{pr-variety} to each $W_{\kappa,i}$ (with $y$
replaced by $(x',y)$) in turn; 
since at most $eK$ such varieties occur, and both $e$ and $K$ are
determined by $\varphi$ only, the total number of Weil numbers that
will occur is bounded in terms of $\varphi$ (for the counting problem)
or $\varphi$ and the degrees of the functions $f$ and $g$ in the
general case.  
\par
The maximal weight is the maximal value of the $2\delta_{\kappa,i}(y)$
of Proposition~\ref{pr-variety} applied to $W_{\kappa,i}$. Denoting it
$2\delta(y)$ it follows that each of the Weil number occurring with
this weight is in fact equal to $q^{\delta(y)}$. 
\par
The remainder of the statement of Theorem~\ref{th-weil-numbers-count}
is clear except maybe that in 
Theorem~\ref{th-weil-numbers-count} the multiplicity
$$
\mu(y)=\triplesum_{\stacksum{\kappa,i,j}
{w(\alpha_{\kappa,i,j}(y))=2\delta(y)}}
{\frac{(-1)^{i+1}}{i!}}
$$
is $>0$ if the set $\varphi(\Fp_q,y)$ is not empty. But otherwise the
number of  elements of $\varphi(\Fp_q,y)$ would be $\ll
q^{\delta(y)-1/2}$, with implied constant depending on $\varphi$ only, 
while, for any $\kappa$, $i$ such that 
$\delta_{\kappa,i}(y)=\delta(y)$, the number of $x\in
W_{\kappa,i,y}(\Fp_q)$ is at least
$\demi q^{\delta(y)}$ if $q$ is large enough (in terms of $\varphi$
only), and the images $\tau_{\kappa,i}(x)$ yield
at least $\tfrac{1}{e}|W_{\kappa,i,y}(\Fp_q)|$ elements of $\varphi(x,y)$
(see~(\ref{eq-imply}) and the property of $e$ following, or look back
at the proof of Theorem~\ref{th-reduction}). Hence we do get the
result stated in Theorem~\ref{th-weil-numbers-count} for all $q$ large
enough. 
\par
There remains to examine the condition stated in
Theorem~\ref{th-weil-numbers-exps} for the exponential sum to have
maximal weight $w(y)<\delta(y)$. For this, assume $w(y)=\delta(y)$. If
$\psi$ is non trivial, then~(\ref{eq-pos-prop}) must hold for some
$\kappa$, $i$ (with $\eta>0$ depending only on $\varphi(x,y)$) with
$\delta_{\kappa,i}(y)=\delta(y)$, i.e., with
$$
|W_{\kappa,i,y}(\Fp_q)|\geq \eta_1 |\varphi(\Fp_q,y)|
$$
where again $\eta_1>0$ depends only on $\varphi(x,y)$. This gives part
(2) of Theorem~\ref{th-weil-numbers-exps}.
\end{proof}

We state finally the following corollary:

\begin{corollary}\label{cor-exp-def}
Let $(f_1,g_1,f_2,g_2,\{\psi\},\{\chi\})$ be the data
defining a family of exponential sums~\emph{(\ref{eq-exp-def})} over
the definable sets $\varphi(\Fp_q,y)$, and let $\gamma(y)$ be as in the
previous statement. There exists $\eta>0$ depending only on
$\varphi(x,y)$ such that for $\psi$ non-trivial we have
$$
S(y,\varphi,\Fp_p)\ll p^{\delta(y)-1/2}
$$
for all primes $p$ and all $y\in\Fp_p^m$ for which there does not
exist $c\in\Fp_p$ such that
$$
\sum_{\stacksum{x\in\varphi(\Fp_p,y)}{f(x)=c}}{1}
\geq \eta \sum_{x\in\varphi(\Fp_p,y)}{1}.
$$
The implied constant depends only on $\varphi(x,y)$ and the degrees of
$f_1$, $f_2$, $g_1$, $g_2$.
\end{corollary}

\begin{proof}
This is immediate after enlarging if necessary the
constant arising from estimating the exponential sum
using~(\ref{eq-decomp}) and $\gamma(y)>0$ for $p$ large enough under
the condition stated.
\end{proof}

\par
The sample statements given in the previous sections
are special
cases of this corollary. In Theorem~\ref{th-one-var}
(with one variable, no parameter, $\psi(x)=e(x/p)$ non trivial), either
$\delta(y)=0$ in  which case~(\ref{eq-one-var}) is trivial, or
$\delta(y)=1$, and if $f$ is a non-constant rational function modulo
$p$, it can not be constant on a set in $\Zz/p\Zz$ containing $\gg p$
points. If $f$ is constant modulo $p$, the estimate is again trivial
as $p\mid N$ in that case.
\par
Theorem~\ref{th-exp-def} on the
other hand is just a rephrasing of the corollary. As for
Theorem~\ref{th-twists}, we are considering the sums
$$
\sum_{x\in\varphi(\Fp_p)}{\chi(g(x))
e\Bigl(\frac{f(x)+\langle h,x\rangle}{p}\Bigr)}
$$
with $\langle h,x\rangle=\sum{h_ix_i}$. Note the degree of $f_h=f+\langle
h,\cdot\rangle$ is uniformly bounded for all $h$. 
Consider an $h$ where the gain for this sum is
$\gamma_h=0$. After
performing the reduction of Theorem~\ref{th-reduction} for
$\varphi(x)$, it is clear from the proof of
Theorem~\ref{th-weil-numbers-exps} above that there is a fixed
$\eta>0$ and a  finite list of subsets $\Phi_i$ of $\varphi(\Fp_p)$
(finite list with cardinality $\leq D$ where $D$ depends only on
$\varphi(x)$) such that if $\gamma_h=0$, then $f$ is constant (say
$=c$) on one of the $\Phi_i$. 
\par
Fixing one $\Phi_i$, which we can assume is not reduced to a
single point (otherwise the desired estimate is trivial), the values
of $h$ corresponding to this $\Phi_i$ are 
such that $\Phi_i$ is contained in an affine hyperplane ``parallel to
$h$'', i.e., with equation $\langle
x,h\rangle=\text{constant}$. Simply by linear algebra, this means that
$\Phi_i$ must be in some intersection $V$ of such affine hyperplanes, and
the assumption that $\Phi_i$ has more than one element means that the
minimal such $V$ is of dimension $k\geq 1$; then one sees that $h$
must lie in a ``complementary'' affine subspace of dimension
$n-k<n$. Hence there will be at most $Dp^{n-1}$ exceptions with
$\gamma_h=0$. (Note that if no $\Phi_i$ is in an hyperplane, we in
fact showed that there are only finitely many ``bad'' values of $h$).

\begin{remark}
The reduction theorem and the subsequent expansions of exponential
sums in terms of Weil numbers suggest an intriguing question: how
intrinsic are those decompositions? Can one define some kind of
cohomology theory for certain sheaves on definable sets over finite
fields in such a way as to obtain formulas like~(\ref{eq-decomp}) as
consequences of a trace formula operating directly at that level?
Since it seems that a condition $q\geq q_0$ is necessary, this may be
better dealt with at the level of pseudo-finite fields. We hope to
come back to such foundational issues.
\par
Finally, the author (at least) can't help wondering if this
introduction of some ideas of logic and model theory might not be one
clue to the hypothetical theory of ``exponential sums over $\Zz$''
that Katz has written about, for instance, in~\cite{katz-esz}.
\par
More down to earth, one may hope that exponential sums estimates for
non-trivial characters could be useful in other areas, noting
for instance that Corollary~\ref{cor-number} of Chatzidakis, van den
Dries and Macintyre has had applications in model
theory (but the author doesn't understand those) and, in remarkable
work of Hrushovski and Pillay~\cite{hru-pillay}, in the study of
algebraic groups over finite fields (in particular the behavior of the
reduction modulo $p$ of a group defined over $\Zz$ for large $p$). 
\end{remark}

\section{Application: definable intervals}
\label{sec-outline}

We come back to general exponential sums~(\ref{eq-gen-exp}), still
with $f(n)=g(n)/p$ for some prime $p$, but now over 
a short interval:
$$
S_{p,\vartheta}(g)=\sum_{0<n\leq N}{e^{2\pi ig(n)/p}}\text{ where }
N=p^{\vartheta}
$$
for some $\vartheta\in ]0,1[$, the inequality $\vartheta<1$ being
characteristic of a ``short'' interval. 
\par
As long as $\vartheta>\tfrac{1}{2}$, a well-known technique of Fourier
completion (see e.g.~\cite[12.2]{ant}) leads to complete sums (i.e.,
with $\vartheta=1$), for which 
the results of algebraic geometry can often by applied, giving a bound
of size roughly $\sqrt{p}(\log p)$, so an estimate for
$S_{p,\vartheta}(g)$ with saving $\theta(N)=N^{1-1/(2\vartheta)}\ra
+\infty$ (see~(\ref{eq-gen-bound})). 
\par
However, very few cases have been handled when $\vartheta\leq
\demi$. Since the most successful high-level approach to exponential
sums in general has been the insertion of a given sum in a sort of
``family'',\footnote{\ Not only for sums over finite fields: the
strategy of the standard methods of Weyl, van der Corput and
Vinogradov can also be understood in this manner.} of sometimes
seemingly unrelated sums, and 
then exploitation of properties of the 
family as a whole to derive individual results, one may hope to do so 
by analogy with the case of complete sums over finite fields by
extending $S_{p,\vartheta}(g)$ in some way to all finite fields
$\Fp_q$.
\par
So the following question is of interest: can one ``lift'' 
the short intervals $0\leq x<p^{\vartheta}$ to subsets of $\Fp_{p^{\nu}}$
in such a way that the corresponding companion sums to $S_p(g)$ have a
sufficiently rich structure to become more accessible?
\par
There are of course many ways to envision such a lifting, but the most
optimistic version is:  are short 
intervals definable by a uniform formula $\varphi(x)$ in one variable
in the language of rings? 
\par
More generally, let $\varphi(x)$ be such a formula. The question is:
when is it the case that
$\varphi(\Fp_p)$ is an \emph{interval} for almost all $p$, i.e.,
when is it the case that 
for all but finitely many $p$ there exist integers $a_p\leq b_p$ for
which 
$$
\varphi(\Fp_p)=\{x\mods{p}\,\mid\, x\equiv i\mods{p}\text{ for some
  integer $i$, $a_p\leq i\leq b_p$}\}.
$$
\par
Of course, given a fixed interval $I=\{n,n+1,\ldots, n+m\}$, $n\geq 0$,
the formula
$$
\varphi(x)\,:\,  (x-n)\cdot (x-(n+1))\cdots (x-(m+n))=0
$$
is such that $\varphi(\Fp_p)$ coincides with the reduction of $I$
modulo $p$ for all $p>m$. So of course similarly the reduction of the
interval $\{n,n+1,\ldots, p-1\}$ of length $p-n$ is defined by
$$
\varphi(x)\,:\, x\cdot (x-1)\cdot (x-(n-1))\not=0.
$$
\par
We can easily use the tools of the previous sections to show that
those are essentially the only possibilities.

\begin{proposition}\label{pr-interval}
Let $\varphi(x)$ be a formula in one variable in the language of
rings. If neither $\varphi(\Fp_p)$ nor $\neg\varphi(\Fp_p)$ are
bounded for all primes, there are only finitely many primes such that 
$\varphi(\Fp_p)$ is the reduction modulo $p$ of an interval.
\end{proposition}

\begin{proof}
By the Main Theorem
of~\cite{zoe} (or Theorem~\ref{th-weil-numbers-count}) applied to
$\varphi$ and $\neg\varphi$, there are
constants $A\geq 1$, $C\geq 0$ and finitely many rationals $\mu_i\in
]0,1[$, such that for each prime $p$ either
$$
|\varphi(\Fp_p)|\leq A\text{ or } |\neg\varphi(\Fp_p)|\leq A
$$
or
$$
\Bigl||\varphi(\Fp_p)|-\mu_i p\Bigr|\leq C\sqrt{p}.
$$
for some $i$. If $p$ runs over a subsequence of primes  (tending to
$+\infty$) for which $\varphi(\Fp_p)$ is an interval, and is
unbounded, the finiteness of the set of $\mu_i$ shows that for some
$i$ a further subsequence exists for which $\varphi(\Fp_p)$ is an
interval of length  $|\varphi(\Fp_p)|\sim \mu_i p$. We will show this
is impossible.
\par
For this we observe first that by Theorem~\ref{th-one-var} we have
\begin{equation}\label{eq-small}
\sum_{x\in \varphi(\Fp_p)}{e\Bigl(\frac{x}{p}\Bigr)}
\ll \sqrt{p},
\end{equation}
for all primes, the implied constant depending only on $\varphi$. 
\par
Let now $p$ be such that $\varphi(\Fp_p)$ is an interval, say
$\varphi(\Fp_p)=\{m_p,\ldots, m_p+n_p-1\}$. Summing a geometric
progression gives 
$$
\Bigl|\sum_{x\in \varphi(\Fp_p)}{e\Bigl(\frac{x}{p}\Bigr)}\Bigl|=
\left|e\Bigl(\frac{m_p}{p}\Bigr)\frac{1-e\Bigl(\dfrac{n_p}{p}\Bigr)}
{1-e\Bigl(\dfrac{1}{p}\Bigr)}
\right|
=\left|\frac{\sin \dfrac{\pi n_p}{p}}{\sin \dfrac{\pi}{p}}\right|.
$$
\par
For the given hypothetical subsequence $p\ra +\infty$,
we have $n_p=|\varphi(\Fp_p)|\sim \mu_i p$,
so the numerator converges 
to $\sin \pi\mu_i\not=0$ (since $\mu_i\notin\{0,1\}$). On the other hand
the denominator is equivalent to $\pi/p$, hence we find that 
\begin{equation}\label{eq-big}
\Bigl|\sum_{x\in \varphi(\Fp_p)}{e\Bigl(\frac{x}{p}\Bigr)}\Bigr|\sim
\pi| \sin \pi\mu_i| p,
\end{equation}
which contradicts~(\ref{eq-small}). This concludes the proof. 
\end{proof}

\begin{remark}
Theorem~\ref{th-one-var} shows more generally that as $p\ra +\infty$
along any sequence with $|\varphi(\Fp_p)|$ unbounded, we have
$$
\lim_{p\ra+\infty}\frac{1}{|\varphi(\Fp_p)|} 
\sum_{x\in \varphi(\Fp_p)}{e\Bigl(\frac{hx}{p}\Bigr)}
=0,
$$
for any fixed integer $h\not=0$.
\par
Hence, using Weyl's equidistribution criterion, it follows that 
the fractional parts $\{\tfrac{x}{p}\}$ for $x\in\varphi(\Fp_p)$ become
equidistributed in $\Rr/\Zz$ for Lebesgue measure as $p$ runs over
primes with $|\varphi(\Fp_p)|\ra +\infty$, i.e., for every
continuous function $f\,:\, \Rr/\Zz\ra \Cc$, one has
$$
\lim_{p\ra +\infty}\frac{1}{|\varphi(\Fp_p)|}
\sum_{x\in\varphi(\Fp_p)}{f\Bigl(\frac{x}{p}\Bigr)}
=\int_{\Rr/\Zz}{f(\theta)d\theta}.
$$
Note that such a result is well-known for
algebraic sets (i.e., formulas without quantifiers); for a much more
general statement of Fouvry and Katz, see~\cite{fouvry-katz}. 
\par
Also, the proposition extends immediately to general arithmetic
progressions : if $\varphi(x)$ is a formula with one variable in the
language of rings, then $\varphi(\Fp_p)$ can be for infinitely many
primes of the form 
$$
A_{a,q,k}=\{a,a+q,\ldots, a+n_pq\}\mods{p}
$$
for some integers $a\geq 1$, $q$ (which may depend on $p$), with
$p\nmid q$, only if either 
$\varphi(\Fp_p)$ or its complement is bounded for all $p$ (the case of
the complement can only occur for $q=1$ of course). Indeed, one need
only compute the exponential sum
$$
\sum_{x\in A_{a,q,k}}{e\Bigl(\frac{\bar{q}x}{p}\Bigr)}
$$
(where $\bar{q}$ is the inverse of $q$ modulo $p$), to obtain an
analogue of~(\ref{eq-big}).
\end{remark}

\section*{Appendix: notation index}

For the reader's convenience, here is a list of the notation
that occur in the reduction theorem and the decomposition theorem for
exponential sums, with a brief explanation of their meaning.
\par
\medskip
\small
\begin{tabular}{lcl}
\notation{$n$}{Number $\geq 0$ of variables in $\varphi(x,y)$}
\notation{$m$}{Number $\geq 0$ of parameters in $\varphi(x,y)$}
\notation{$K$}{Theorem~\ref{th-reduction} (integer $\geq 1$, number of
disjunctions expressing $\varphi$)}
\notation{$s$}{Theorem~\ref{th-reduction} (integer $\geq 0$, number
of auxiliary parameters needed)}
\notation{$x'$}{Value or variables for the auxiliary parameters}
\notation{$\kappa$}{Index running from $1$ to $K$}
\notation{$\Phi_{\kappa}$}{Formulas in the disjunction expressing
$\varphi$}
\notation{$k$}{Proof of Theorem~\ref{th-reduction} (integer $\geq 0$,
number of existential terms in $\Phi_{\kappa}$)}
\notation{$r$}{Proof of Theorem~\ref{th-reduction} (integer $\geq 0$,
number of terms in formulas $\Phi_{\kappa}$)}
\notation{$f_{\kappa,\cdot}$}{Terms occurring in
  $\Phi_{\kappa}$} 
\notation{$h_{\kappa,\cdot}$}{Terms existing in
  existential form in $\Phi_{\kappa}$}
\notation{$\Psi_{\kappa}$}{Proof of Theorem~\ref{th-reduction}
  (auxiliary quantifier free formulas)} 
\notation{$q_0,p_0$}{Generically, value of $q$ or $p$ so that a
  statement holds for $q\geq q_0$ or $p\geq p_0$}
\notation{$e$}{Theorem~\ref{th-reduction} (integer $\geq 1$, maximal
  number of pre-images for a given $x\in \varphi(\Fp_q,y)$)}
\notation{$i,j$}{Indices running from $1$ to $e$}
\notation{$W_{\kappa,i}$}{Affine schemes projecting ``to $\varphi$''}
\notation{$\pi_{\kappa,i}$}{Projection from $W_{\kappa,i}$ to space of
  parameters}
\notation{$\tau_{\kappa,i}$}{Projection from $W_{\kappa,i}$ to
  definable sets}
\notation{$\eps(\alpha)$}{Paragraph before
  Proposition~\ref{pr-variety} (``sign'' of a Weil number)}
\notation{$B_1$}{Theorem~\ref{th-weil-numbers-count} (maximal number
of Weil numbers for the counting function)}
\notation{$B_2$}{Theorem~\ref{th-weil-numbers-exps} (maximal number of
 Weil numbers for exponential sum)}
\notation{$\beta$}{Theorems~\ref{th-weil-numbers-count}
and~\ref{th-weil-numbers-exps} (number of Weil numbers occurring
for given parameter $y$)}
\notation{$j$}{Theorems~\ref{th-weil-numbers-count}
and~\ref{th-weil-numbers-exps} (index running from $1$ to $\beta$)}
\notation{$\alpha_{\kappa,i,j}(y)$}{Theorems~\ref{th-weil-numbers-count}
  and~\ref{th-weil-numbers-exps} (Weil numbers giving exponential
  sum)}
\notation{$\delta(y)$}{Theorem~\ref{th-weil-numbers-count} (maximal
  weight of Weil numbers in counting points)}
\notation{$\mu(y)$}{Theorem~\ref{th-weil-numbers-count} (density of
  point counting)}
\notation{$w(y)$}{Theorem~\ref{th-weil-numbers-exps} (maximal weight
  of Weil numbers in exponential sum)}
\notation{$\gamma(y)$}{Theorem~\ref{th-weil-numbers-exps} (gain in
  bound for exponential sum)}
\end{tabular}

\end{document}